\documentclass[11pt,a4paper,leqno]{article}
\usepackage{amsfonts,amssymb,mathrsfs,amsmath}

\usepackage{amscd}

\vspace{0.5cm}

 4

\font\erm=cmr8

\usepackage{color}
\usepackage{graphicx}

\author{Andrzej Krzysztof Kwa\'sniewski}
\title{Report  on a new type  - mixed  V and U binomials' recurrence.}

\newtheorem{defn}{Definition}

\newtheorem{theoremn}{Theorem}

\newcommand{\fnomialF}[3]{ {{#1} \choose {#2}}_{\!\!#3} }

\begin{document}

\begin{center}
\noindent {\sffamily Report  on a new type  - mixed  V and U binomials' recurrence.}  \\ 

\vspace{0.5cm} 
\sffamily Andrzej Krzysztof Kwa\'sniewski

\vspace{0.5cm}
{\erm
\sffamily Member of the Institute of Combinatorics and its Applications, 
Winnipeg, Manitoba, Canada \\
	PL-15-674 Bia\l ystok, Konwaliowa 11/11, Poland\\
	e-mail: kwandr@gmail.com\\}
\end{center}

\vspace{0,5cm}

\noindent \textbf{Summary} 

\noindent \sffamily In \cite[2011]{2011 akk} we had proposed \textit{H(x)}-binomials' recurrence formula appointed by Ward-Horadam $ H(x)= \left\langle H_n(x)\right\rangle _{n\geq 0}$ functions' sequence i.e. any functions' sequence solution of the second order recurrence with functions' coefficients.\\ 

\noindent Then we had delivered  a new type \textit{H(x)}-binomials' "mixed" recurrence for primordial case.
\noindent Here  we report  on this new type  - mixed  $V$ and $U$ binomials' recurrence - in brief. 
\vspace{0,1cm}

\noindent For the abundant list of references on the generalized binomials recurrences subject   we refer the reader to \cite[2011]{2011 akk}.

\vspace{0,2cm}
\noindent AMS Classification Numbers: 05A10 , 05A30.

\vspace{0,1cm}

\vspace{0.4cm}


\section{Sketchy introduction to the history of the subject}

\vspace{0,1cm}

\noindent Up to our knowledge it was Fran\c{c}ois \'Edouard Anatole Lucas in \cite[1878]{EdL} who was the first who had  \textit{not only} defined \textit{fibonomial} coefficients as stated in \cite[1989]{K-W}   by Donald  Ervin  Knuth and Herbert Saul Wilf but who was the first who had defined  \textit{$U_n \equiv n_{p,q}$}-binomial coefficients ${n \choose k}_U \equiv {n \choose k}_{p,q}$  and had derived a recurrence for them: see page 27, formula (58) \cite[1878]{EdL}. Then - referring to Lucas - the investigation relative to  divisibility properties of relevant number Lucas sequences $D$, $S$ as well as numbers' $D$ - binomials and  numbers' $D$ - multinomials was continued in \cite[1913]{Carmichel} by  Robert Daniel Carmichel;  see pp. 30,35 and 40 in \cite[1913]{Carmichel}  for $ U \equiv D = \left\langle D_n\right\rangle_{n\geq 0}$ and ${n \choose {k_1,k_2,...,k_s}}_D$ - respectively. Note there also formulas (10), (11) and (13) which might perhaps  serve to derive explicit untangled form of recurrence for the  $ \textit{V}$- \textit{binomial coefficients} ${n \choose k}_V \equiv {n \choose k}_S$ denoted by  primordial Lucas sequence $\left\langle S_n\right\rangle_{n\geq 0}= S \equiv V$. 
\vspace{0,1cm}

\noindent Multinomial coefficients' \textit{recurrences} are not present in that early and other works and up  to our knowledge a special case of such appeared at first in \cite[1979]{Shannon 1979 multi} by  Anthony G. Shannon. More on that - in what follows after Definition 3?.

\vspace{0,1cm}


\noindent \textbf{1.2. We deliver here}   - a new  \textit{H(x)}-binomials' recurrence formula appointed by Ward-Horadam $ H(x)= \left\langle H_n(x)\right\rangle _{n\geq 0}$ field of zero characteristic nonzero valued functions' sequence which comprises  for $H \equiv H(x=1)$ number sequences case - the  \textit{V}-binomials' recurrence formula determined by the primordial Lucas sequence of the second kind  $V = \left\langle V_n\right\rangle_{n\geq 0}$  \cite[2011]{2011 akk} as well  as its well elaborated companion fundamental Lucas sequence of the first kind $ U = \left\langle U_n\right\rangle_{n\geq 0}$ which gives rise in its turn to the  \textit{U}-binomials' recurrence as in \cite[1878]{EdL} , \cite[1949]{JM},  \cite[1964]{TF},  \cite[1969]{Gould},  \cite[1989]{G-V}  or in  \cite[1989]{K-W}  and so on. 

\vspace{0,1cm}

\noindent We do it by following recent applicable work \cite[2009]{Savage}  by Nicolas A. Loehr and Carla  D.  Savage (see more in \cite[2011]{2011 akk}).

\vspace{0,1cm}


\noindent \textbf{This looked for} here new \textit{H(x)}-binomials' overall recurrence formulas -(\emph{recall}: these should be encompassing \textit{V}-binomials for primordial Lucas sequence $V$) - this looked for here recurrence for $V$ is not present neither in \cite{EdL} nor in \cite{Savage}, nor in \cite[1915]{Fon} , nor in \cite[1936]{Ward}, nor in \cite[1949]{JM}. Neither we find it in  all 206 quoted in \cite[2011]{2011 akk} references.\\

\vspace{0,1cm}

\noindent \textbf{\textcolor{red}{Interrogation}}. Might it be so that such an  overall recurrence formula including  $V$ does not exist? The method to obtain formulas is more general than actually  produced family of formulas?  See Theorems 2a and 2b, and Example 1  in what follows.


\section{Prerequisites}

\vspace{0.1cm}

\noindent The Lucas sequence   $V = \left\langle V_n\right\rangle_{n\geq 0}$  is called the Lucas sequence of the second kind - see: \cite[1977, Part I]{KKK}, or \textbf{primordial} - see
\cite[1979]{ShaHor}.\\     

\noindent The Lucas sequence  $ U = \left\langle U_n\right\rangle_{n\geq 0}$   is called the Lucas sequence of the first kind - see: \cite[1977, Part I]{KKK}, or \textbf{fundamental}  - see p. 38 in \cite[1949]{JM} or see \cite[1979]{Shannon 1979 multi}  and  \cite[1979]{ShaHor}.     

\vspace{0.1cm} 

\noindent In the sequel - following  \cite[2011]{2011 akk} - we shall deliver the looked for recurrences for $H$-binomial coefficients ${n \choose k}_H$  determined by the Ward-Horadam sequence $H$ - defined below. It appears to be  mixed  $V$ and $U$ binomials' recurrence in the primordial case  $H = V$.

\vspace{0.1cm}

\noindent In compliance with Edouard Lucas' \cite[1878]{EdL} and twenty, twenty first century $p,q$-people's notation \cite[2011]{2011 akk} we  at first review here in brief the general second order recurrence;
(compare this review  with the  recent "Ward-Horadam" peoples' paper  \cite [2009]{He-Siue} by Tian-Xiao He and Peter Jau-Shyong Shiue or earlier  $p,q$-papers \cite[2001]{Sun-Hu-Liu} by Zhi-Wei Sun, Hong Hu, J.-X. Liu  and \cite[2001]{Hu-Sun}  by Hong Hu and Zhi-Wei Sun). And with respect to natation:  If in  \cite[1878]{EdL} Fran\c{c}ois  \'Edouard Anatole Lucas  had been used $a=\textbf{p}$ and $b = \textbf{q}$ notation,  he would be perhaps at first glance notified and recognized as a Great Grandfather of all the $(p,q)$ - people. Let us start then introducing reconciling and matched denotations and nomenclature.

\begin{equation}
H_{n+2} = P \cdot H_{n+1}  - Q \cdot H_n ,\ \;n\geq 0 \ and \;  H_0 = a,\; H_1 = b.                                                                                     
\end{equation}
which is sometimes being written  in $\left\langle P,-Q \right\rangle \mapsto  \left\langle s,t \right\rangle$ notation.

\begin{equation}
H_{n+2} = s \cdot H_{n+1}  + t\cdot H_n ,\ \;n\geq 0 \ and \;  H_0 = a,\; H_1 = b.                                                                                     
\end{equation}
We exclude the cases when (2) is recurrence of the first order, therefore we assume that the roots $p,q$ of (5) are distinct   $p\neq q$  and $\frac{p}{q}$ is not the root of unity. We shall come back to this finally while formulating this note  observation named Theorem 2a.

\vspace{0.2cm}
\noindent Simultaneously and collaterally  we  mnemonically pre adjust the starting point to discuss the $F(x)$ polynomials' case via - if entitled - antecedent "$\mapsto$ action": 
$H \mapsto H(x)$, $s \mapsto s(x)$, $t \mapsto t(x)$, etc. 

\begin{equation}
H_{n+2}(x) = s(x) \cdot H_{n+1}(x) + t(x)\cdot H_n, \; n\geq 0, \ H_0 = a(x), H_1 = b(x).                                                                                     
\end{equation}

\vspace{0.2cm}
\noindent enabling recovering  explicit formulas also for sequences of  polynomials correspondingly generated by the above linear recurrence of order 2  - with Pafnuty Lvovich  Tchebysheff  polynomials and the generalized Gegenbauer-Humbert polynomials included. See for example Proposition 2.7 in the  recent "Ward-Horadam peoples'" paper  \cite [2009]{He-Siue}  by  Tian-Xiao He and Peter Jau-Shyong Shiue.

\vspace{0.2cm}

\noindent The general solution of (1):   $ H(a,b;P,Q) = \left\langle H_n\right\rangle _{n\geq 0}$ is being called  throughout  this note - \textbf{Ward-Horadam number'sequence} \cite[2011]{2011 akk}.

\vspace{0.2cm}

\noindent The Lucas $U$-binomial coefficients ${n \choose k}_U \equiv {n \choose k}_{p,q}$  are then defined as follows: (\cite[1878]{EdL},  \cite[1915]{Fon},  \cite[1936]{Ward},  \cite[1949]{JM},  \cite[1964]{TF}, \cite[1969]{Gould} etc.)

\vspace{0.3cm}

\begin{defn}
	Let  $U$ be as in \cite[1878]{EdL} i.e $U_n \equiv n_{p,q}$ then $U$-binomial coefficients for any $n,k \in \mathbb{N}\cup\{0\}$ are defined as follows
	\begin{equation}
		{n \choose k}_U \equiv {n \choose k}_{p,q} = \frac{n_{p,q}!}{k_{p,q}! \cdot (n-k)_{p,q}!} = \frac{n_{p,q}^{\underline{k}}}{k_{p,q}!}
	\end{equation}
	\noindent where $n_{p,q}! = n_{p,q}\cdot(n-1)_{p,q}\cdot ... \cdot 1_{p,q}$ and $n_{p,q}^{\underline{k}} = n_{p,q}\cdot(n-1)_{p,q}\cdot ...\cdot (n-k+1)_{p,q}$  and ${n \choose k}_U = 0$ for $k > n$.
\end{defn} 

\vspace{0.3cm}

\begin{defn}
Let $V$ be as in \cite[1878]{EdL} i.e $V_n = p^n +q^n$, hence $V_0 = 2$ and  $V_1 = p + q = s$. Then $V$-binomial coefficients for any $n,k \in \mathbb{N}\cup\{0\}$ are defined as follows
	\begin{equation}
		{n \choose k}_V =\frac{V_n!}{V_k!\cdot V_(n-k)!} = \frac{V_n^{\underline{k}}}{V_k!}
	\end{equation}
	\noindent where $V_n! = V_n \cdot V_{n-1}\cdot...\cdot V_1$ and $V_n^{\underline{k}}=V_n \cdot V_{n-1}\cdot ...\cdot V_{n-k+1}$ and ${n \choose k}_V = 0$ for $k > n$.
\end{defn}

\vspace{0.2cm}

\noindent One automatically generalizes number $F$-binomial coefficients' array to functions $F(x)$-\textbf{multinomial} coefficients' array (see  \cite[2004]{akk3 2004} and references to umbral calculus therein) while for \textit{number sequences} $F = F(x=1)$ the $F$-multinomial coefficients see p. 40 in \cite[1913]{Carmichel} by  Robert Daniel Carmichel , see  \cite[1936]{Ward} by Morgan Ward and \cite[1969]{Gould}  by Henri W. Gould; (see more in \cite[2011]{2011 akk}).

\vspace{0.2cm}

\noindent \textbf{$\psi(x)$-multinomial coefficients.}

\vspace{0.1cm}

\noindent Considerations in \cite[1979]{Shannon 1979 multi}  by  Anthony G. Shannon  and an application in \cite[2001]{Rich 2001} by Thomas M. Richardson as well as relevance to  umbral calculus 
\cite[2003]{akk ITSF 2003} or  \cite[2001]{akk ITSF 2001} constitute  motivating circumstances for considering now \textit{functions'} $F(x)$-\textbf{binomial}  and  $F(x)$-\textbf{multinomial} coefficients. This is the case ($F(x)=\psi(x)$) for example in  \cite[2004]{akk3 2004}  \textbf{wherein we read:} 

$$
\left(x_{1}+_{\psi}x_{2}+_{\psi} \ldots +_{\psi}x_{k}\right)^{n}=
\sum_{\begin{array}{l} s_{1},\ldots s_{k}=0\\
s_{1}+s_{2}+\ldots +s_{k}=n
\end{array}}^{n}\binom{n}{s_{1},\ldots ,s_{k}}_{\psi}x_{1}^{s_{1}}\ldots x_{k}^{s_{k}}
$$
where
$$\binom{n}{s_{1},\ldots ,s_{k}}_{\psi}=\frac{n_{\psi}!}{(s_{1})_{\psi}!\ldots (s_{k})_{\psi}!}.$$
Here above   the $u$-{\em multinomial} number sequence formula from \cite[1936]{Ward} by Morgan Ward is extended mnemonically to $\psi(x)$ function sequence definition of shifting $x$ arguments formula written in Kwa\'sniewski upside-down notation (see for example \cite[2003]{akk ITSF 2003},  \cite[2001]{akk ITSF 2001}, \cite[2002]{akk Lodz2002},  \cite[2005]{akk Allahabad  2005} for more ad this notation).    

\noindent  The above formulas introduced  in case of number sequences in \cite[1936]{Ward} by Morgan Ward were recalled from \cite[1936]{Ward}  by Alwyn F. Horadam and  Anthony G.  Shannon  in \cite[1976]{HorSha}  were the Ward Calculus of sequences framework (including umbral derivative and corresponding exponent) was used to enunciate two types of  Ward's Staudt-Clausen theorems pertinent to this general calculus of number sequences.

\section{$H(x)$-binomial and mixed binomial coefficients' recurrence}

\vspace{0.2cm}

\noindent \textbf{3.1.}  Let us recall convention resulting from (3).

\vspace{0.2cm}
\noindent Recall. The general solution of (3):   $$ H(x) \equiv H(a(x),b(x);s(x),t(x)) = \left\langle H_n(x)\right\rangle _{n\geq 0}$$ 
is being called  throughout this paper - \textbf{Ward-Horadam functions' sequence}.

\vspace{0.2cm}

\begin{theoremn}
Let us admit shortly the abbreviations: $g_k(r,s)(x) = g_k(r,s)$ , $k=1,2$.  Let $s,r>0$. Let $F(x)$ be any zero characteristic field nonzero valued functions' sequence ($F_n(x) \neq 0$). Then

\begin{equation}
{r+s \choose r,s}_{F(x)}=g_1(r,s)\cdot{r+s-1 \choose r-1,s}_{F(x)}+g_2(r,s)\cdot{r+s-1 \choose r,s-1}_{F(x)}
\end{equation}
where   $ {r \choose r,0}_{F(x)} = {s \choose 0,s}_{F(x)} =1$ and 
\begin{equation}
 F(x)_{r+s} =  g_1(r,s) \cdot F(x)_r   +  g_2(r,s) \cdot F(x)_s.                               
\end{equation}
are equivalent.
\end{theoremn}

\vspace{0.3cm}

\noindent \textbf{An on the way historical note}
\noindent  Donald  Ervin  Knuth   and    Herbert Saul Wilf in \cite[1989]{K-W} stated that  Fibonomial coefficients and the recurrent relations for them appeared already in 1878  Lucas work (see: formula (58)  in \cite[1878]{EdL}  p. 27 ;  for $U$-binomials which "Fibonomials" are special case of). More over on this very p. 27 Lucas formulated a conclusion from his (58) formula which may be  stated in notation of this paper formula (2) as follows: \textit{if}  $s,t  \in\mathbb{Z}$ and $H_0 =0$ , $H_1 = 1$  \textit{then} $H \equiv U$  and $ \fnomialF{n}{k}{U} \equiv \fnomialF{n}{k}{n_{p,q}} \in\mathbb{Z}$.  

\vspace{0.1cm}

\noindent Consult for that also the next century references: \cite[1910]{Bachmann 1910}  by  Paul Gustav Heinrich Bachmann  or later on - \cite[1913]{Carmichel}  by  Robert Daniel Carmichel  [p. 40]  or  \cite[1949]{JM} by Dov Jarden and Theodor Motzkin where in all quoted positions it was also shown that $n_{p,q}$ - binomial coefficients are integers  - for $p$ and $q$ representing  distinct roots of the characteristic equation for (3) - see below.

\vspace{0.2cm}

\noindent It seems to be the right place now to underline that the  \textit{addition formulas} for Lucas sequences below with respective hyperbolic trigonometry formulas and also consequently $U$-binomials'recurrence formulas - stem from commutative  ring $R$ identity: $(x-y)\cdot (x+y) \equiv x^2 - y^2, x,y \in R$.

\vspace{0.3cm}

\noindent Indeed. Recall the characteristic equation notation:  $\textcolor{blue}{z^2 = s \cdot z + t}$ as is common in many publications with the restriction:  $p,q$ roots are distinct. Recall that $p+q=s$  and $p\cdot q=-t$. Let $\Delta = s^2 - 4t$. Then $\Delta = (p-q)^2$. Hence we have as in \cite[1878]{EdL}

\begin{equation}
2\cdot U_{r+s} =  U_r V_s  +  U_s V_r ,\ \ \ \ 
2\cdot V_{r+s} =  V_r V_s  +  \Delta \cdot U_s U_r. 
\end{equation}
Taking here into account the  $\textbf{U}$-\textcolor{blue}{addition formula} i.e. the first of two trigonometric-like $L$-addition formulas (42) from \cite[1878]{EdL} ($L = U,V$)  
one readily recognizes that  the $U$-binomial recurrence from the Corollary 18 in  \cite[2009]{Savage} is  the $U$-binomial recurrence (58) \cite[1878]{EdL}  which may be rewritten after Fran\c{c}ois  \'Edouard Anatole Lucas in multinomial notation and stated as follows: \textit{according  to the Theorem 2a  below the following is true}:

$$2\cdot U_{r+s} =  U_r V_s  +  U_s V_r$$ 
\textit{is equivalent to}

\begin{equation}
2\cdot{r+s \choose r,s}_{n_{p,q}}= V_s \cdot{r+s-1 \choose r-1,s}_{n_{p,q}}+V_r \cdot{r+s-1 \choose r,s-1}_{n_{p,q}}.
\end{equation}

\vspace{0.2cm}

\noindent However there is no companion  $V$-binomial recurrence i.e. for ${r+s \choose r,s}_V$ neither in  \cite[1878]{EdL}  nor in \cite[2009]{Savage} as well as all other quoted papers - up to knowledge of this note author.  

\vspace{0.1cm}

\noindent \textbf{The End} of \textit{the on the way historical note}.

\vspace{0.4cm}

\noindent The looked for  $H(x)$-binomial recurrence  (6) accompanied by (7) might be  then given right now in the form of (10) adapted to  - Ward-Lucas functions'sequence case notation while keeping in mind that of course the expressions for $h_k(r,s)(x)$,  $k=1,2$ below are designated by this $ F(x) = H(x)$ choice and as a matter of fact are appointed by the recurrence (3). 

\vspace{0.3cm}

\noindent For the sake of commodity to write down an observation to be next let us admit shortly the abbreviations: $h_k(r,s)(x) = h_k(r,s)= h_k$ , $k=1,2$. Then for $H(x)$ of the form (21?) we evidently have what follows.

\vspace{0.6cm}

\noindent \textbf{\textcolor{blue}{Theorem 2a}.}

\begin{equation}
{r+s \choose r,s}_{H(x)} = h_1(r,s){r+s-1 \choose r-1,s}_{H(x)} + h_2(r,s){r+s-1 \choose r,s-1}_{H(x)},
\end{equation}
where  ${r \choose r,0}_{H(x)} = {s \choose 0,s}_{H(x)} =1,$  \textit{is equivalent to}

\begin{equation}
 H_{r+s}(x) =  h_1(r,s)H_r(x)   +  h_2(r,s) H_s(x).                               
\end{equation}
where  $H_n(x)$ is explicitly given by (21)  and (22)  formulas in \cite[2011]{2011 akk}.\\ 
\vspace{0.1cm}
\noindent \textbf{The end} of the Theorem 2a.

\vspace{0.4cm}


\noindent Let us now come back to consider the case of $V$-binomials recurrence.  
For that to do carefully let us at first make precise the main item. $ H(x) = \left\langle H_n(x)\right\rangle_{n\geq 0}$ is considered as a solution of \textcolor{red}{second} order recurrence (3) with peculiar case of (3) becoming the first order recurrence \textcolor{blue}{excluded}. This is equivalent to say that  $0\neq p \neq q \neq 0 $ and  $A \neq 0 \neq B$     (compare with \cite[1974]{Hilton 1974 partition}  by Anthony J. W. Hilton), where for notation convenience we shall again use awhile shortcuts  for  (3):  

$$a(x) \equiv a , b(x) \equiv b , s(x) \equiv s , t(x) \equiv t ,$$

$$p(x)\cdot q(x) = - t(x) \equiv p \cdot q = - t,$$

$$ H(x)= H(A(x),B(x),p(x),q(x)) = H(A,B) =A \cdot H(1,1;) ,$$
$$ H(1,1) = V(x)=V = (p-q)U(1,-1) , \  U = U(x)= U(1,1) , U(A,B) = U(A,B,p,q)$$
referring to this note (23) and the next to (23) abbreviations: 
$$H(x) \equiv H(a(x),b(x);s(x),t(x)) \equiv  H(A(x),B(x);s(x),t(x)),$$   $$H_n(x) \equiv H_n(a(x),b(x);s(x),t(x)),$$  
$$H_n(x) = H_n(A,B) = A\cdot p^n + B\cdot q^n \ H_n(1,1) = V_n(x)= V_n.$$

\vspace{0.2cm}

\noindent Now although $H(A,B) =(p-q)\cdot U(A,-B)$  the corresponding recurrences are different. Recall and then compare corresponding recurrences:

\noindent \textit{The identity }  
$$(p-q)\cdot(p^{r+s} - q^{r+s}) \equiv (p^{s+1} - q^{s+1})\cdot(p^r - q^r) - p\cdot q (p^{r-1} - q^{r-1}) \cdot(p^s - q^s)$$
due to $p\cdot q = -t$  is \textit{equivalent to }

$$
U_{r + s}(p,q) = U_{s+1} \cdot U_r(p,q) + t \cdot U_{r-1}\cdot U_s(p,q)
$$
(combinatorial derivation - see \cite[1999]{1BQ 1999},\cite[1999]{2BQ 1999}). This recurrence in its turn \noindent\textit{is equivalent to} 

$$
{r+s \choose r,s}_U = U_{s+1} \cdot {r+s-1 \choose r-1,s}_U + t \cdot U_{r-1}\cdot{r+s-1 \choose r,s-1}_U.
$$
\noindent Similarly \textit{the identity }  
$$(p-q)\cdot(p^{r+s} + q^{r+s}) \equiv (p^{s+1} - q^{s+1})\cdot(p^r + q^r) - p\cdot q (p^{r-1} + q^{r-1}) \cdot(p^s - q^s)$$
due to $p\cdot q = -t$  is \textit{equivalent to}

\begin{equation}
V_{r + s}(p,q) = U_{s+1} \cdot V_r(p,q) + t \cdot V_{r-1}\cdot U_s(p,q)
\end{equation}
For combinatorial interpretation derivation of the above for Fibonacci and Lucas sequences ($t=1=s$) see \cite[1999]{1BQ 1999}, \cite[1999]{2BQ 1999} and then see for more \cite[2003]{BQ2003} also by Arthur T. Benjamin and Jennifer J. Quinn. A fundamental progress in combinatorial interpretation of generalized binomials  was made in \cite[2009]{BSCS} by Bruce E. Sagan and Carla D. Savage.

\begin{defn}
Let  $ \left\{ n(x) \right\} \equiv U_n $ . $\left\langle n(x) \right\rangle \equiv V_n = p^n(x) +q^n(x)$, hence $V_0 = 2$ and  $V_1 = p + q = s(x)$. Let  $(p(x)- q(x)) \cdot U_n = p^n(x) - q^n(x)$, hence $U_0 = 0$ and  $V_1 = 1$; (roots are distinct). Then $V$-mixed-$U$  binomial coefficients for any $r,s \in \mathbb{N} \cup \{0\}$ ; are defined as follows
\begin{equation}
		{r+s \choose r,s}_{\left\langle . \right\rangle / \left\{ . \right\}} = \frac{V_{r+s}!}{V_r!\cdot U_s!}, 
\end{equation}
\end{defn}
${n \choose k}_{\left\langle . \right\rangle / \left\{ . \right\}}= 0 $  for $k>n$  and  ${n \choose 0}_{\left\langle . \right\rangle/\left\{ . \right\}}= 1.$

\vspace{0.2cm}

\noindent Note that:  ${r+s \choose r,s}_{\left\langle .\right\rangle/\left\{.\right\}} \neq {r+s \choose s,r}_{\left\langle .\right\rangle/\left\{.\right\}}.$

\vspace{0.2cm}

\noindent Note that: ${r+s \choose r,s}_{\left\langle .\right\rangle/\left\{.\right\}} \neq {r+s \choose r,s}_{ \left\{.\right\}/ \left\langle .\right\rangle }$

\vspace{0.3cm}

\noindent \textbf{\textcolor{red}{Theorem 2b}.} The recurrence (12) \noindent\textit{is equivalent to} 

\begin{equation}
{r+s \choose r,s}_{\left\langle . \right\rangle / \left\{ . \right\}} = U_{s+1} \cdot {r+s -1  \choose r-1,s}_{\left\langle . \right\rangle / \left\{ . \right\}} + t \cdot U_s \cdot {r+s-1 \choose r,s-1}_{\left\langle . \right\rangle / \left\{ . \right\}}.
\end{equation}

\vspace{0.3cm}

\noindent \textbf{Example 1.} Application of  Theorems 1, 2a and  the Theorem 2b method.  
\noindent Recall the characteristic equation notation:  $z^2 =\textcolor{blue}{z^2 = s \cdot z + t}$  with the restriction that $p,q$ roots are distinct. Recall that $p+q=s$  and $p\cdot q=-t$. Let $\Delta = s^2 - 4t$. Then $\Delta = (p-q)^2$. Hence we have as in \cite[1878]{EdL}

$$
2\cdot U_{r+s} =  U_r V_s  +  U_s V_r ,\ \ \ \ 
2\cdot V_{r+s} =  V_r V_s  +  \Delta \cdot U_s U_r. 
$$
Consider   $\textbf{V}$-\textcolor{blue}{addition formula} i.e. the second  of two trigonometric-like addition formulas (42) from \cite[1878]{EdL}. One readily recognizes that  \textit{according  to the Theorem 2b  the following is true}:     $2 \cdot V_{r+s} =  V_r V_s  +  \Delta \cdot U_s U_r$  \textit{is equivalent to}

\begin{equation}
2\cdot {r+s \choose r,s}_{\left\langle . \right\rangle / \left\{ . \right\}} = V_s \cdot {r+s -1  \choose r-1,s}_{\left\langle . \right\rangle / \left\{ . \right\}} + \Delta \cdot U_r \cdot {r+s-1 \choose r,s-1}_{\left\langle . \right\rangle / \left\{ . \right\}}.
\end{equation}
The effort to discover combinatorial interpretation of mixed binomials is  part of a dream of a forthcoming "in statu nascendi"  note.\\
\noindent s for mixed $F(x)$-\textbf{multinomial} coefficients' arrays - there are many kinds of them. This is also is a part of a subject  of a forthcoming "in statu nascendi"  note.

\vspace{0.3cm}


\end{document}